\documentstyle[times,aps,amsthm,amsmath,amsfonts,%
    enumerate,twocolumn,floats]{revtex}

\newcommand{\cN}{{\cal N}}
\newcommand{\cO}{{\cal O}}
\newcommand{\cP}{{\cal P}}
\newcommand{\cQ}{{\cal Q}}

\newcommand{\all}{{\mbox{All}}}
\newcommand{\none}{{\mbox{none}}}
\newcommand{\tP}{{\bf 0}}
\newcommand{\tN}{{\bf 1}}
\newcommand{\tO}{{\bf 2}}
\newcommand{\tM}{{\bf 3}}
\newcommand{\tQ}{{\mbox{\boldmath $\infty$}}}
\newcommand{\near}{{2_1}}
\newcommand{\far}{{2}}

\newtheorem{claim}{Claim}
\theoremstyle{definition}
\newtheorem{question}{Question}

\begin{document}

\title{Three-player impartial games}
\date{November 10, 1998}
\author{James Propp}
\address{Department of Mathematics, University of Wisconsin}
\maketitle
\begin{abstract}
Past efforts to classify impartial three-player combinatorial games
(the theories of Li~\cite{Li} and Straffin~\cite{Str})
have made various restrictive assumptions about 
the rationality of one's opponents
and the formation and behavior of coalitions.
One may instead adopt an agnostic attitude towards such issues,
and seek only to understand in what circumstances
one player has a winning strategy
against the combined forces of the other two.
By limiting ourselves to this more modest theoretical objective,
and by regarding two games as being equivalent
if they are interchangeable in all disjunctive sums
as far as single-player winnability is concerned,
we can obtain an interesting analogue of Grundy values
for three-player impartial games.
\end{abstract}

\section{Introduction}

Let us begin with a very specific problem:
Assume $G$ is an impartial (positional) game
played by three people
who alternate moves in cyclic fashion
(Natalie, Oliver, Percival, Natalie, Oliver, Percival, ...),
under the convention that
the player who makes the last move wins.
Let $H$ be another such game.
Suppose that the second player, Oliver,
has a winning strategy for $G$.
Suppose also that Oliver
has a winning strategy for $H$.
Is it possible for Oliver
to have a winning strategy
for the disjunctive sum $G+H$ as well?

Recall that an {\bf impartial positional game}
is specified by
(i) an initial position,
(ii) the set of all positions that can arise during play, and
(iii) the set of all legal moves from one position to another.
The winner is the last player to make a move.
To avoid the possibility of a game going on forever,
we require that from no position may there be
an infinite chain of legal moves.
The disjunctive sum of two such games $G,H$
is the game in which a legal move
consists of making a move in $G$ (leaving $H$ alone)
{\it or\/} making a move in $H$ (leaving $G$ alone).
Readers unfamiliar with the theory of two-player impartial games
should consult \cite{BCG} or \cite{Con}.

It is important to notice that in a three-player game,
it is possible that {\it none\/} of the players has a winning strategy.
The simplest example is the Nim game that starts from the position $1+2$,
where 1 and 2 denote Nim-heaps of size one and two respectively.
As usual, a legal move consists of taking
a number of counters from a single heap.
In this example, the first player has no winning move,
but his actions determine whether
the second or third player will win the game.
None of the players has a winning strategy.
That is, any two players can cooperate
to prevent the remaining player from winning.
It is in a player's interest
to join such a coalition of size two
if he can count on his partner
to share the prize with him --
unless the third player counters
by offering an even bigger share of the prize.
This kind of situation is well known in
the theory of ``economic'' (as opposed to positional) games.
In such games, however, play is usually simultaneous rather than sequential.

Bob Li \cite{Li} has worked out a theory of multi-player positional games
by decreeing that a player's winnings depend on how recently
he has moved when the game ends
(the last player to move wins the most,
the player who moved before him wins the next most,
and so on),
and by assuming that each player will play rationally
so as to get the highest winnings possible.
Li's theory, when applied to games like Nim,
leads to quite pretty results,
and this is perhaps sufficient justification for it;
but it is worth pointing out
that, to the extent that game theory is supposed to be applicable
to the actual playing of games,
it is a bit odd to assume
that one's adversaries are going to play perfectly.
Indeed, the only kind of adversaries a sensible person would play with,
at least when money is involved,
are those who do {\it not\/} know the winning strategy.
Only in the case of two-player games
is it the case that a player has a winning strategy 
against an arbitrary adversary
if and only if he has a winning strategy
against a perfectly rational adversary.

Phil Straffin \cite{Str} has his own approach
to three-player games.
He adopts a policy (``McCarthy's revenge rule'')
governing how a player should act
in a situation where he himself cannot win
but where he can choose
which of his opponents will win.
Straffin analyzes Nim under such a revenge rule,
and his results are satisfying if taken on their own terms,
but the approach is open to the same practical objections as Li's.
Specifically,
if a player's winning strategy depends on the assumption
that his adversaries will be able to recognize when they can't win,
then the player's strategy is guaranteed to work
only when 
his opponents can see all the way to the leaves of the game tree.
In this case, at least one of them (and perhaps each of them)
believes he can't win;
so why is he playing?

The proper response to such objections,
from the point of view of someone who wishes to understand 
real-world games,
is that theories like Li's and Straffin's
are prototypes of more sophisticated theories,
not yet developed,
that take into account the fact that players of
real-life games
are partly rational and partly emotional creatures,
capable of such things as stupidity and duplicity.

It would be good to have a framework into which
the theories of Li and Straffin,
along with three-player game-theories of the future,
can be fitted.
This neutral framework would make no special assumptions
about how the players behave.
Here, we develop such a theory.
It is a theory designed to answer the single question
``Can I win?,''
asked by a single player
playing against two adversaries of unknown characteristics.
Not surprisingly,
the typical answer given by the theory is ``No'';
in most positions,
any two players can gang up on the third.
But it turns out that there is a great deal to be said
about those games in which one of the players
{\it does\/} have a winning strategy.

In addition to the coarse classification of three-player games
according to who (if anyone) has the winning strategy,
one can also carry out a fine classification of games
analogous to, but much messier than,
the classification of two-player games according to Grundy-value.
The beginnings of such a classification
permit one to answer the riddle with which this article opened;
the later stages lead to many interesting complications
which have so far resisted all attempts at comprehensive analysis.

\section{Notation and Preliminaries}

Games will be denoted by the capital letters $G$, $H$, $X$, and $Y$.
As in the two-player theory, we can assume that
every position carries along with it the rules of play to be applied,
so that each game may be identified with its initial position.
The game $G'$ is an {\bf option} of $G$
if it is legal to move from $G$ to $G'$.
To build up all the finite games,
we start from the {\bf null-game} $0$
(the unique game with no options)
and recursively define
$G = \{G'_1, G'_2, ...\}$
as the game with options $G'_1, G'_2, ...$.
The game $\{0\}$ will be denoted by 1,
the game $\{0,1\}$ will be denoted by 2, and so on.
(It should always be clear from context
whether a given numeral
denotes a number or a Nim game.)

We recursively define the relation of {\bf identity}
by the rule that $G$ and $H$ are identical if and only if
for every option $G'$ of $G$
there exists an option $H'$ of $H$
identical to it, and vice versa.
We define {\bf (disjunctive) addition},
represented by $+$,
by the rule that $G+H$ is the game
whose options are precisely the games of the form
$G'+H$ and $G+H'$.
It is easy to show that identity is an equivalence relation
that respects the ``bracketing'' and addition operations,
that addition is associative and commutative,
and that $0$ is an additive identity.

The following abbreviations will prove convenient:
$$\begin{array}{rll}
GH & \mbox{means} & G+H \\
G^n & \mbox{means} & G+G+...+G \ \ \mbox{($n$ times)} \\
m_n & \mbox{means} & \{\{\cdots\{m\}\cdots\}\} \ \ \mbox{($n$ layers deep)}
\end{array}$$

Thus, $\{12\}^3 4_5$ denotes
$$\{1+2\} + \{1+2\} + \{1+2\} + \{\{\{\{\{4\}\}\}\}\}.$$
(We'll never need to talk about Nim-heaps of size $>9$,
so our juxtaposition convention won't cause trouble.)
Note that for all $G$, the games
$G0$, $G^1$, $G_0$, and $G$ are identical.

Relative to any non-initial position in the course of play,
one of the players has just moved (the Previous player)
and one is about to move (the Next player);
the remaining player is the Other player.
At the start of the game,
players Next, Other, and Previous
correspond to the first, second, and third players
(even though, strictly speaking,
there was no ``previous'' move).
We call $G$ a Next-game ($\cN$-game)
if there is a winning strategy for Next,
and we let $\cN$ be the set of $\cN$-games;
$\cN$ is the {\bf type} of $G$,
and $G$ {\bf belongs} to $\cN$.
We define $\cO$-games and $\cP$-games in a similar way.

If none of the players has a winning strategy,
we say that $G$ is a Queer game ($\cQ$-game).
In a slight abuse of notation,
I will often use ``$=$'' to mean ``belongs to'',
and use the letters $N,O,P,Q$
to stand for unknown games belonging to
these respective types.
Thus I will write $1=N$, $11=O$, $111=P$, etc.;
and the problem posed in the Introduction
can be formulated succinctly as:
solve $O+O=O$ or prove that no solution exists.
(At this point I invite the reader to tackle $Q+Q=O$.
There is a simple and elegant solution.)

The following four rules provide a recursive method
for classifying a game:

\begin{itemize}
\item[(1)]
$G$ is an $\cN$-game exactly if it has some $\cP$-game as an option.
\item[(2)]
$G$ is an $\cO$-game exactly if all of its options are $\cN$-games,
{\it and} it has at least one option (this proviso prevents us from
mistakenly classifying $0$ as an $\cO$-game).
\item[(3)]
$G$ is a $\cP$-game exactly if all of its options are $\cO$-games.
\item[(4)]
$G$ is a $\cQ$-game exactly if none of the above conditions is satisfied.
\end{itemize}

Using these rules, it is possible to analyze a game completely
by classifying all the positions in its game-tree, from leaves to root.  

\section{Some Sample Games}

Let us first establish the types of the simpler Nim games.
It's easy to see that
\begin{align*}
0   & = P, \\
1   & = \{0\} = \{P\} = N, \\
11  & = \{1\} = \{N\} = O, \\
111 & = \{11\} = \{O\} = P,
\end{align*}
and so on;
in general, the type of $1^n$ is $\cP$, $\cN$, or $\cO$
according as the residue of $n$ mod 3 is 0, 1, or 2.
Also
\begin{align*}
2 & = N, \\
3 & = N,
\end{align*}
and so on, because in each case Next can win by taking the whole heap.
\begin{align*}
12    & = \{1,2,11\} = \{N,N,O\} = Q, \\
112   & = \{11,12,111\} = \{O,Q,P\} = N, \\
1112  & = \{111,112,1111\} = \{P,N,N\} = N, \\
11112 & = \{1111,1112,11111\} = \{N,N,O\} = Q,
\end{align*}
and so on; in general, the type of $1^n 2$ is $\cN$, $\cQ$, or $\cN$
according as the residue of $n$ mod 3 is 0, 1, or 2.
The winning strategy for these $\cN$-games is simple:
reduce the game to one of the $\cP$-positions $1^{3k}$.

$1+1=11$ is a solution of the equation $N+N=O$.
Does $G=N$ imply that $G+G=O$ in general?
We can easily see that the answer is ``No'':
$$ 2+2 = 22 = \{12,2\} = \{Q,N\} = Q. $$
($12$ is identical to $21$, so they can be treated as a single option.)

Here are some more calculations which will be useful later.
\begin{align*}
\{2\}    & = \{N\} = O \\
\{\{2\}\}& = \{O\} = P \\
\{1,11\} & = \{N,O\} = Q \\
\{2,11\} & = \{N,O\} = Q
\end{align*}

\section{Adding Games}
\label{s:adding}

The type of $G+H$ is not in general determined by
the types of $G$ and $H$.
(For example, 1 and 2 are both of type $\cN$,
but $1+1=O$ while $2+2=Q$.)
That is,
addition does not respect the relation ``belongs to the same type as''.
To remedy this situation we define {\bf equivalence} ($\equiv$) by
the condition that $G \equiv H$ if and only if
for all games $X$, $G+X$ and $H+X$ belong to the same type.
It is easy to show that ``equivalence'' is an equivalence relation,
that it respects bracketing and addition,
and that if $G' \equiv H'$ then $\{G',H',...\} \equiv \{H',...\}$
(that is, equivalence options of a game may be conflated).

We are now in a position to undertake the main task of this section:
determining the addition table.
Recall that in the two-player theory, there are only two types
($\cN$ and $\cP$) and their addition table is as shown in
Table \ref{t:two}.

\begin{figure}[htb]
\begin{center}
\begin{tabular}{c|ll}
$+$ & $P$ & $N$ \\ \hline 
$P$ & $P$ & $N$ \\
$N$ & $N$ & $PN$
\end{tabular}
\end{center}
\caption{The two-player addition table.}
\label{t:two}
\end{figure}

\noindent
Here, the entry $PN$ denotes the fact
that the sum of two $\cN$-games 
can be either a $\cP$-game
or an $\cN$-game.

The analogous addition table for three-player games
is given by Table \ref{t:three}.

\begin{figure}[htb]
\begin{center}
\begin{tabular}{c|llll}
$+$ & $P$  & $N$   & $O$   & $Q$   \\ \hline
$P$ & $PQ$ & $NQ$  & $OQ$  & $Q$   \\
$N$ & $NQ$ & $NOQ$ & $PNQ$ & $NQ$  \\
$O$ & $OQ$ & $PNQ$ & $NQ$  & $NOQ$ \\
$Q$ & $Q$  & $NQ$  & $NOQ$ & $OQ$   
\end{tabular}
\end{center}
\caption{The three-player addition table.}
\label{t:three}
\end{figure}

Notice that in one particular case (namely $G=P$ and $H=Q$, or vice versa),
knowing the types of $G$ and $H$ does tell one which type $G+H$ belongs to,
namely $Q$.
A corollary of this is that $P+P+...+P+Q=Q$.

To prove that Table~\ref{t:two} applies,
one simply finds solutions of the allowed ``equations'' 
$P+P=P$, $P+N=N$ (from which $N+P=N$ follows), $N+N=P$, and $N+N=N$,
and proves that the forbidden equations $P+P=N$ and $P+N=P$ have no solutions.
To demonstrate the validity of Table \ref{t:three},
we must find solutions to twenty-two such equations,
and prove that the remaining eighteen have no solutions.

Table~\ref{t:sums} shows the twenty-two satisfiable equations 
and their solutions.

\begin{figure}[htb]
\begin{center}
\begin{tabular}{@{\extracolsep{3pt}}ll|ll}
Equation & Solution & Equation & Solution \\ \hline
$P+P=P$ & $0+0$   & $P+P=Q$ & $\{\{2\}\}+\{\{2\}\}$ \\
$N+P=N$ & $1+0$   & $N+P=Q$ & $1111+\{\{2\}\}$      \\
$N+N=O$ & $1+1$   & $O+P=Q$ & $\{2\}+111$           \\
$O+P=O$ & $11+0$  & $O+N=Q$ & $1+\{2\}$             \\
$O+N=P$ & $11+1$  & $O+O=Q$ & $\{2\}+\{2\}$         \\
$O+O=N$ & $11+11$ & $N+N=N$ & $112+1$               \\
$Q+P=Q$ & $12+0$  & $N+O=N$ & $2+11$                \\
$Q+N=Q$ & $12+2$  & $Q+N=N$ & $12+1$                \\
$Q+O=Q$ & $22+11$ & $Q+O=N$ & $12+11$               \\
$Q+Q=Q$ & $12+12$ & $Q+O=O$ & $\{2,11\} + 11$       \\
$N+N=Q$ & $2+2$   & $Q+Q=O$ & $\{1,11\} + \{1,11\}$
\end{tabular}
\end{center}
\caption{Some sums.}
\label{t:sums}
\end{figure}

And now, the proofs of impossibility
for the eighteen impossible cases.

\begin{claim} None of the following is possible.
\begin{align}
O+P&=N \label{opn} \\
N+P&=P \label{npp} \\
O+O&=P \label{oop} \\
P+P&=O \label{ppo} \\
O+N&=O \label{ono}
\end{align}
\end{claim}

\begin{proof} By (joint) infinite descent.
Here, as in subsequent proofs,
the infinite-descent ``boilerplating'' is omitted.

Note that none of the hypothetical $\cP$-games
in equations \eqref{opn}-\eqref{ppo} can be the $0$-game,
so all of these games $X,Y$ have options.

Suppose \eqref{opn} holds; say $X=O$, $Y=P$, $X+Y=N$.
Some option $X'+Y$ or $X+Y'$ must be a $\cP$-game.
But then we have either 
$N+P=P$ (every option $X'$ must be an $\cN$-game), which is \eqref{npp}, 
or $O+O=P$ (every option $Y'$ must be an $\cO$-game), which is \eqref{oop}.

Suppose \eqref{npp} holds; say $X=N$, $Y=P$, $X+Y=P$.
Then there exists $X'=P$, which must satisfy $X'+Y=P+P=O$ (equation \eqref{ppo}).

Suppose \eqref{oop} holds; say $X=O$, $Y=O$, $X+Y=P$.
Then there exists $Y'=N$, which must satisfy $X+Y'=O+N=O$ (equation \eqref{ono}).

Suppose \eqref{ppo} holds; say $X=P$, $Y=P$, $X+Y=O$.
Then there exists $X'=O$, which must satisfy $X'+Y=O+P=N$ (equation \eqref{opn}).

Finally, suppose \eqref{ono} holds; say $X=O$, $Y=N$, $X+Y=O$.
Then there exists $Y'=P$, which must satisfy $X+Y'=O+P=N$ (equation \eqref{opn}).
\end{proof}

\begin{claim} None of the following is possible.
\begin{align}
P+P&=N \label{ppn} \\
O+P&=P \label{opp} \\
N+P&=O \label{npo}
\end{align}
\end{claim}

\begin{proof} By infinite descent. A solution to \eqref{ppn} yields an
(earlier-created) solution to \eqref{opp}, which yields a solution to
\eqref{npo}, which yields a solution to \eqref{ppn}.
\end{proof}

\begin{claim} It is impossible that
\begin{equation} N+N=P \label{nnp} \end{equation}
\end{claim}

\begin{proof} By contradiction.  A solution to \eqref{nnp} would
yield a solution to \eqref{npo}.
\end{proof}

\begin{claim} None of the following is possible.
\begin{align}
Q+P&=N \label{qpn} \\
Q+P&=P \label{qpp} \\
Q+O&=P \label{qop} \\
Q+P&=O \label{qpo} \\
Q+N&=O \label{qno}
\end{align}
\end{claim}

\begin{proof} By infinite descent (making use of earlier results as well).

Suppose \eqref{qpn} holds with $X,Y$.
Some option $X'+Y$ or $X+Y'$ must be a $\cP$-game.
In the former event, we have
$X' \neq P$ (since $X=Q$),
so that either $N+P=P$ (equation \eqref{npp}),
$O+P=P$ (equation \eqref{opp}), or
$Q+P=P$ (equation \eqref{qpp}); 
in the latter event we have $Q+O=P$ (equation \eqref{qop}).

Suppose \eqref{qpp} holds with $X,Y$.
Since $X=Q$, it has an option $X'$ of type $\cN$ or type $\cQ$
(for if all options of $X$ were $\cO$-games and $\cP$-games,
$X$ would be of type $\cP$ or $\cN$).
If $X'=N$, then we have $X'+Y=N+P=O$ (equation \eqref{npo}),
and if $X'=Q$, then we have $X'+Y=Q+P=O$ (equation \eqref{qpo}).

Suppose \eqref{qop} holds with $X,Y$.
Then $X+Y'=Q+N=O$ (equation \eqref{qno}).

Suppose \eqref{qpo} holds with $X,Y$.
Since $X=Q$, it has an option $X'$ of type $\cO$ or of type $\cQ$
(for if all options of $X$ were $\cN$-games and $\cP$-games,
$X$ would be of type $\cO$ or $\cN$).
$X'=O$ yields $X'+Y=O+P=N$ (equation \eqref{opn}),
and $X'=Q$ yields $X'+Y=Q+P=N$ (equation \eqref{qpn}).

Finally, suppose \eqref{qno} holds with $X,Y$.
Then there exists $Y'=P$,
which must satisfy $X+Y'=Q+P=N$ (equation \eqref{qpn}).
\end{proof}

\begin{claim} It is impossible that
\begin{equation} Q+N=P \label{qnp} \end{equation}
\end{claim}

\begin{proof} By contradiction.  A solution to \eqref{qnp} would yield
a solution to \eqref{qpo}.
\end{proof}

\begin{claim} Neither of the following is possible:
\begin{align}
Q+Q&=N \label{qqn} \\
Q+Q&=P \label{qqp}
\end{align}
\end{claim}

\begin{proof} By infinite descent.

Suppose \eqref{qqn} holds with $X,Y$.
Then some option of $X+Y$ must be a $\cP$-game;
without loss of generality, we assume $X+Y'=P$.
But $X=Q$, and we have already ruled out
$Q+P=P$ (equation \eqref{qpp}),
$Q+N=P$ (equation \eqref{ono}),
and $Q+O=P$ (equation \eqref{qop}),
so we have $X+Y'=Q+Q=P$ (equation \eqref{qqp}).

Suppose \eqref{qqp} holds with $X,Y$.
$X$ must have an $\cN$-option or $\cQ$-option $X'$,
but if $X'=N$ then $X'+Y=N+Q=O$ (equation \eqref{qno}),
which can't happen; so $X'=Q$.
Similarly, $Y$ has a $\cQ$-option $Y'$.
$X'+Y=O$, so $X'+Y'=Q+Q=N$ (equation \eqref{qqn}).

(Note that the second half of this proof
requires us to look two moves ahead,
rather than just one move ahead as in the preceding proofs.)
\end{proof}

The remaining case is surprisingly hard to dispose of;
the proof requires us to look {\it five\/} moves ahead.

\begin{claim} It is impossible that
\begin{equation} O+O=O \label{ooo} \end{equation}
\label{c:impossible}
\end{claim}

\begin{proof} By infinite descent.
Suppose \eqref{ooo} holds with $X,Y$.
For all $X'$ we have $X'+Y=N$,
so that $X'+Y$ must have some $\cP$-option;
but this $\cP$-option cannot be of the form $X'+Y'$,
since $N+N \neq P$ (equation \eqref{nnp}).
Hence there must exist an option $X''$ of $X'$
such that $X''+Y=P$.
This implies that $X''=N$, since none of the cases 
$O+O=P$ (equation \eqref{oop}), $P+O=P$ (equation \eqref{opp}), $Q+O=P$
(equation \eqref{qop}) 
can occur.
Similarly, every $Y'$ has an option $Y''$ such that
$X+Y''=P$, $Y''=N$.
Since $X''+Y$ is a $\cP$-game,
$X''+Y'$ and $X'+Y''$ are $\cO$-games
and $X''+Y''$ is an $\cN$-game.
One of the options of $X''+Y''$ must be a $\cP$-game;
without loss of generality, say $X'''+Y''=P$.
Since $Y''=N$ and since none of the cases 
$N+N=P$ (equation \eqref{nnp}), $P+N=P$ (equation \eqref{npp}),
$Q+N=P$ (equation \eqref{qnp}) 
can occur,
$X'''$ must be an $\cO$-game.
But recall that $X''+Y$ is a $\cP$-game,
so that its option $X'''+Y$ is an $\cO$-game.
This gives us $X'''+Y=O+O=O$,
which is an earlier-created solution to \eqref{ooo}.
\end{proof}

The proof of Claim \ref{c:impossible} completes the proof of the validity of Table~\ref{t:three}.
Observe that this final clinching claim,
which answers the article's opening riddle in the negative,
depends on five of the preceding six claims.
Our straightforward question thus seems to lack a straightforward solution.
In particular, one would like to know of
a winning strategy for the Natalie-and-Percival coalition in the game $G+H$
that makes use of Oliver's winning strategies for $G$ and $H$.
Indeed, it would be desirable to have strategic ways of understanding
all the facts in this section.

At this point it is a good idea to switch to a notation
that is more mnemonically helpful than $N$, $O$, and $P$,
vis-\`a-vis addition.
Let $\tP$, $\tN$, and $\tO$
denote the Nim-positions $0$, $1$, $11$, respectively.
Also, let $\tQ$ be the Nim-position $22$.
(Actually, we'll want these symbols to represent
the equivalence classes of these respective games,
but that distinction is unimportant right now.)
We will say that two games $G$, $H$ are {\bf similar}
if they have the same type;
in symbols, $G \sim H$.
Every game is thus similar to exactly one
of $\tP$, $\tN$, $\tO$, and $\tQ$.
We can thus use these four symbols to classify our games by type;
for instance, instead of writing $G=N$,
we can write $G \sim \tN$.

Here is the rule for
recursively determining the type of a game
in terms of the types of its options,
restated in the new notation:
\begin{itemize}
\item[(1)]
$G$ is of type $\tN$ exactly if it has some option of type $\tP$.
\item[(2)]
$G$ is of type $\tO$ exactly if all of its options are of type $\tN$,
{\it and} it has at least one option.
\item[(3)]
$G$ is of type $\tP$ exactly if all of its options are of type $\tO$.
\item[(4)]
$G$ is of type $\tQ$ exactly if none of the above conditions is satisfied.
\end{itemize}

Here is the new addition table for 3-player game types;
it resembles a faulty version of the modulo 3 addition table.

\begin{figure}[htb]
\begin{center}
\begin{tabular}{c|cccc}
$+$   & $\tP$    & $\tN$       & $\tO$       & $\tQ$       \\  \hline
$\tP$ & $\tP\tQ$ & $\tN\tQ$    & $\tO\tQ$    & $\tQ$       \\
$\tN$ & $\tN\tQ$ & $\tN\tO\tQ$ & $\tP\tN\tQ$ & $\tN\tQ$    \\
$\tO$ & $\tO\tQ$ & $\tP\tN\tQ$ & $\tN\tQ$    & $\tN\tO\tQ$ \\
$\tQ$ & $\tQ$    & $\tN\tQ$    & $\tN\tO\tQ$ & $\tO\tQ$   
\end{tabular}
\end{center}
\caption{The new, improved three-player addition table.}
\label{t:improved}
\end{figure}

It is also worthwhile to present the ``subtraction table''
as an object of study in its own right.
To this end define $\tM = 111$ as an alternative to $\tP$.

\begin{figure}[htb]
\begin{center}
\begin{tabular}{c|cccc}
$-$   & $\tP$  & $\tN$  & $\tO$       & $\tQ$     \\ \hline
$\tM$ & $\tM$  & $\tO$  & $\tN$       & $\none$   \\
$\tO$ & $\tO$  & $\tN$  & $\tP\tQ$    & $\tO \tQ$ \\
$\tN$ & $\tN$  & $\all$ & $\tN\tO\tQ$ & $\tN \tO$ \\
$\tQ$ & $\all$ & $\all$ & $\all$      & $\all$
\end{tabular}
\end{center}
\caption{The three-player ``subtraction'' table.}
\label{t:subtraction}
\end{figure}

The minuend is indicated by the row
and the subtrahend by the column.
Note that subtraction is not a true operation on games;
rather, the assertion ``$\tN - \tO$ is $\tN \tO \tQ$''
means that if $G,H$ are games
such that $G+H \sim \tN$ and $G \sim \tO$
then $H \sim$ $\tN$, $\tO$, or $\tQ$.

The six entries in the upper left corner of the subtraction table
(the only entries that are single types)
correspond to assertions that can be proved by joint induction
without any reference to earlier tables.
In fact, a good alternative way to prove 
that addition satisfies Table \ref{t:improved}
would be to prove that addition satisfies 
the properties implied by the six upper-left entries
in Table \ref{t:subtraction}
(by joint induction)
and then to prove three extra claims:
(i)   if $G \sim \tO$ and $H \sim \tO$ then $G+H \not\sim \tO$;
(ii)  if $G \sim \tQ$ and $H \sim \tQ$ then $G+H \not\sim \tP$; and
(iii) if $G \sim \tQ$ and $H \sim \tQ$ then $G+H \not\sim \tN$.

\section{Adding Games to Themselves}

Another sort of question related to addition
concerns the disjunctive sum of a game with itself.
Recall that in two-player game theory,
a strategy-stealing argument can be used to show
that the sum of a game of type $\cN$ with itself
must be of type $\cP$
(even though a sum of two distinct games of type $\cN$
can be of either type $\cP$ or type $\cN$).

We seek a similar understanding of
what happens when we add a three-player game to itself.
Table \ref{t:doubling} shows the possible types $G+G$ can have
in our three-player theory,
given the type of $G$.

\begin{figure}[htb]
\begin{center}
\begin{tabular}{c|c}
$G$   & $G+G$    \\ \hline
$\tP$ & $\tP\tQ$ \\
$\tN$ & $\tO\tQ$ \\
$\tO$ & $\tN\tQ$ \\
$\tQ$ & $\tO\tQ$
\end{tabular}
\end{center}
\caption{The doubling table.}
\label{t:doubling}
\end{figure}

To verify that all the possibilities listed here can occur,
one can simply look at the examples given
at the beginning of Section \ref{s:adding}.
To verify that none of the omitted possibilities can occur,
it {\it almost\/} suffices to consult Table \ref{t:improved}.
The only possibility that is not ruled out by the addition table
is that there might be a game $X$ with
$X \sim \tN$, $X+X \sim \tN$.

Suppose $X$ were such a game.
Then $X$ would have to have a $\cP$-option $X'_1$
(now we call it a $\tP$-option)
along with another option $X'_2$ such that $X+X'_2 \sim \tP$.
This implies that $X'_1+X'_2 \sim \tO$ and $X'_2+X'_2 \sim \tO$.
Since $X'_1 \sim \tP$,
the condition $X'_1+X'_2 \sim \tO$
implies (by way of Table \ref{t:improved})
that $X'_2 \sim \tO$.
But $X'_2+X'_2 \sim \tO$
implies (by way of Table \ref{t:improved})
that $X'_2 \sim \tN \mbox{ or} \tQ$.
This contradiction shows that no such game $X$ exists,
and completes the verification of Table \ref{t:doubling}.

In the same spirit, we present a trebling table (Table \ref{t:trebling}),
showing the possible types $G+G+G$ can have
given the type of $G$.

\begin{figure}[htb]
\begin{center}
\begin{tabular}{c|c}
$G$   & $G+G+G$  \\ \hline
$\tP$ & $\tP\tQ$ \\
$\tN$ & $\tP\tQ$ \\
$\tO$ & $\tP\tQ$ \\
$\tQ$ & $\tQ$
\end{tabular}
\end{center}
\caption{The trebling table.}
\label{t:trebling}
\end{figure}

To prove that all the possibilities listed in the first three rows
can actually occur,
one need only check that 
$0 + 0 + 0 \sim \tP$,
$\{\{2\}\} + \{\{2\}\} + \{\{2\}\} \sim \tQ$,
$1 + 1 + 1 \sim \tP$,
$2 + 2 + 2 \sim \tQ$,
$11 + 11 + 11 \sim \tP$,
and $\{2\} + \{2\} + \{2\} \sim \tQ$.
To prove that the nine cases not listed cannot occur
takes more work.

Four of the cases are eliminated
by the observation that
$G+G+G$ can never be of type $\tN$
(the second and third players can always make the Next player lose
by using a copy-cat strategy).
Tables 3 and 5 allow one to eliminate three more cases.
The next two claims take care of the final two cases.

\begin{claim} If $G \sim \tQ$, then $G+G+G \not\sim \tO$.
\end{claim}

\begin{proof} 
Suppose $X\sim\tQ$ with $X+X+X\sim\tO$.
Let $X^\alpha$ be an option of $X$.
Since $X^\alpha+X+X\sim\tN$,
$X^\alpha+X+X$ must have a $\tP$-option
of the form $X^\alpha+X^\beta+X$
(for $X^\beta$ some option of $X$)
or of the form $X^{\alpha\gamma}+X+X$
(for $X^{\alpha\gamma}$ some option of $X^\alpha$).
In either case, we find that the $\tQ$-game $X$,
when added to some other game
($X^\alpha+X^\beta$ or $X^{\alpha\gamma}+X$),
yields a game of type $\tP$;
this is impossible, by Table \ref{t:improved}.
\end{proof}

\begin{claim} If $G \sim \tO$, then $G+G+G \not\sim \tO$.
\end{claim}

\begin{proof}
Suppose $X\sim\tO$ with $X+X+X\sim\tO$.
Notice that $X'+X+X\sim\tN$
for every option $X'$ of $X$.

Case I:
There exist options $X^\alpha$, $X^\beta$ of $X$
(possibly the same option)
for which $X^\alpha + X^\beta + X \sim \tP$.
Then its option $X^\alpha + (X^\beta + X^\beta) \sim \tO$.
Since $X^\alpha \sim \tN$, Table \ref{t:subtraction} gives
$X^\beta + X^\beta \sim \tN$.
But this contradicts Table \ref{t:doubling}, since $X^\beta \sim \tN$.

Case II:
There do not exist two such options of $X$.
Let $X^\alpha$ be an option of $X$.
Since $X^\alpha+X+X \sim \tN$,
and since there exists no $X^\beta$
for which $X^\alpha + X^\beta + X \sim \tP$,
there must exist an option $X^{\alpha\gamma}$ of $X^\alpha$
such that $X^{\alpha\gamma}+X+X \sim \tP$.
$X+X \sim \tN \mbox{ or } \tQ$,
by Table \ref{t:doubling},
but $X+X$ cannot be of type $\tQ$,
since adding $X^{\alpha\gamma}$
yields a $\tP$-position.
Hence $X+X \sim \tN$,
and Table \ref{t:subtraction} implies $X^{\alpha\gamma} \sim \tO$.
Since $X+X \sim \tN$,
there must exist an option $X^\delta$
with $X^\delta + X \sim \tP$.
Everything we've proved so far about $X^\alpha$
applies equally well to $X^\delta$
(since all we assumed about $X^\alpha$
was that it be some option of $X$).
In particular, $X^\delta$
must have an option $X^{\delta\epsilon}$
such that
$X^{\delta\epsilon} \sim \tO$.
However, since
$X^{\delta\epsilon}+X$ is an option of
the $\tP$-position $X^\delta + X$,
$X^{\delta\epsilon}+X \sim \tO$.
Hence $X^{\delta\epsilon}$ and $X$
are two $\tO$-positions
whose sum is a $\tO$-position,
contradicting Table \ref{t:improved}.
\end{proof}

\section{Nim for Three}

We wish to classify all Nim-positions
as belonging to $\cN$, $\cO$, $\cP$, or $\cQ$ ---
or rather, as we now put it,
as being similar to $\tP$, $\tN$, $\tO$, or $\tQ$.
We will actually do more,
and determine the equivalence classes of Nim games.
Table \ref{t:basic} shows the games we have classified so far
(on the left) and their respective types (on the right).

\begin{figure}[htb]
\begin{center}
$\begin{array}{lcl}
0 & \ \ & (\tP) \\
1, 11, 111, 1111, ... & \ \ & (\tN,\tO,\tP,\tN,\tO,\tP,...) \\
2, 12, 112, 1112, ... & \ \ & (\tN,\tQ,\tN,\tN,\tQ,\tN,...) \\
3 & \ \ & (\tN) \\
22 & \ \ & (\tQ)
\end{array}$
\end{center}
\caption{Basic positions of Nim.}
\label{t:basic}
\end{figure}

\noindent
We will soon see that every Nim-game
is equivalent to one of the Nim-games in Table \ref{t:basic}.
We call these {\it reduced\/} Nim-positions.
The last paragraph of this section gives a procedure
for converting a three-player Nim-position
into its reduced form.

Throughout this section (and the rest of this article),
the reader should keep in mind the difference between
the notations 2 and $\tO$.  
The former is a single Nim-heap of size 2;
the latter is the game-type that corresponds to a second-player win.
Note in particular that 2 is not of type $\tO$
but rather of type $\tN$.

We start our proof of the validity of Table \ref{t:basic} by showing that
no two games in the table are equivalent to each other.
In this we will be assisted by Tables \ref{t:type1} and \ref{t:type2}.

\begin{figure}[htb]
\begin{center}
$\begin{array}{r|cccccccccccc}
m\backslash n &0 &1 &2 &3 &4 &5 &6 &7 &8 &9 &10 &\ldots \\ \hline
0 &\tN &\tO &\tP &\tN &\tO &\tP &\tN &\tO &\tP &\tN &\tO &\ldots \\
1 &\tQ &\tQ &\tN &\tO &\tP &\tN &\tO &\tP &\tN &\tO &\tP &\ldots \\
2 &\tN &\tQ &\tQ &\tQ &\tN &\tO &\tP &\tN &\tO &\tP &\tN &\ldots \\
3 &\tN &\tQ &\tQ &\tQ &\tQ &\tQ &\tN &\tO &\tP &\tN &\tO &\ldots \\
4 &\tQ &\tQ &\tQ &\tQ &\tQ &\tQ &\tQ &\tQ &\tN &\tO &\tP &\ldots \\
5 &\tN &\tQ &\tQ &\tQ &\tQ &\tQ &\tQ &\tQ &\tQ &\tQ &\tN &\ldots \\
\vdots & \vdots & & & & & \vdots & & & & & \vdots & \\
\end{array}$
\end{center}
\caption{The type of $1^m + 2_n$.}
\label{t:type1}
\end{figure}

Table \ref{t:type1} gives the types for games of the form $1^m+2_n$.
Each row of the chart gives what we shall call
the {\it signature\/} of $1^m$, relative to
the sequence $2, \{2\}, \{\{2\}\}, ...\:$.
Since no two games of the form $1^m$ have the same signature,
no two are equivalent.

Similarly, Table \ref{t:type2} is
the signature table for games of the form $1^m 2$,
relative to $2_n$.

\begin{figure}[htb]
\begin{center}
$\begin{array}{r|cccccccccccc}
m\backslash n& 0 &1 &2 &3 &4 &5 &6 &7 &8 &9 &10 &\ldots \\\hline
0 &\tQ &\tQ &\tN &\tQ &\tN &\tN &\tQ &\tN &\tN &\tQ &\tN &\ldots \\
1 &\tQ &\tQ &\tQ &\tQ &\tN &\tQ &\tN &\tN &\tQ &\tN &\tN &\ldots \\
2 &\tQ &\tQ &\tQ &\tQ &\tQ &\tQ &\tN &\tQ &\tN &\tN &\tQ &\ldots \\
3 &\tQ &\tQ &\tQ &\tQ &\tQ &\tQ &\tQ &\tQ &\tN &\tQ &\tN &\ldots \\
4 &\tQ &\tQ &\tQ &\tQ &\tQ &\tQ &\tQ &\tQ &\tQ &\tQ &\tN &\ldots \\
5 &\tQ &\tQ &\tQ &\tQ &\tQ &\tQ &\tQ &\tQ &\tQ &\tQ &\tQ &\ldots \\
\vdots &\vdots & & & & &\vdots & & & & \vdots& & \\
\end{array}$
\end{center}
\caption{The type of $1^m 2 + 2_n$.}
\label{t:type2}
\end{figure}

We see that all the games $1^m$ and $1^m 2$ are distinct.

What about $22$?
It can't be equivalent to $1^{3k+1} 2$ for any $k$
(even though both are $\tQ$-games),
because $22+1 \sim \tQ$
while $1^{3k+1} 2 + 1 = 1^{3k+2} 2 \sim \tN$.

What about 3?
It can't be equivalent to $1^{3k+1}$ for any $k$,
because $3+1 \sim \tQ$ while $1^{3k+1} + 1 \sim \tO$;
it can't be equivalent to $1^{3k} 2$
because $3+2_2 \sim \tN$
while $1^{3k} 2 + 2_2 \sim \tQ$;
it can't be equivalent to $1^{3k+2} 2$
because $3+1 \sim \tQ$ 
while $1^{3k+2} 2 + 1 = 1^{3k+3} 2 \sim \tN$;
and it can't be equivalent to $2$
because $\{0,11\} + 2 \sim \tO$ 
while $\{0,11\}+3 \sim \tQ$.

Now that we know that all of the Nim games 
in Table \ref{t:basic} are inequivalent,
let us show that {\it every\/} Nim game
is equivalent to one of these.

\begin{claim} $mn \sim \tQ$ for all $m,n \geq 2$.
\label{c:gangup}
\end{claim}

\begin{proof} Any two players can gang up on the third,
by depleting neither heap
until the victim has made his move,
and then removing both heaps.
\end{proof}

\begin{claim} The following are true for all games $G$:
\begin{enumerate}[(a)]
\item $Gn \not\sim \tP$ for $n \geq 2$.
\item $Gn \not\sim \tO$ for $n \geq 3$.
\item If $Gm \sim \tN$ then $Gn \sim \tN$, for $m,n \geq 2$.
\item $G1n \not\sim \tO$ for $n \geq 2$.
\item $Gmn \not\sim \tN$ for $m,n \geq 2$.
\item $Gmn \not\sim \tO$ for $m,n \geq 2$.
\item $Gmn \not\sim \tP$ for $m,n \geq 2$.
\end{enumerate}
\label{c:inequiv}
\end{claim}

\begin{proof}\ \par

(a) Suppose $Gn\sim\tP$.
Then its options $G1$ and $G$ are $\tO$-games.
But since $G$ is also an option of $G1$,
this is a contradiction.

(b) Suppose $Gn\sim\tO$.
Then $G$, $G1$, and $G2$
are all $\tN$-games,
and in particular $G2$ must have a $\tP$-option.
That $\tP$-option can be neither $G$ nor $G1$,
so there must exist $G'2\sim\tP$,
contradicting (a).

(c) Assume $Gm\sim\tN$.
Then either $G\sim\tP$ or $G1\sim\tP$
(no other option of $Gm$ can be of type $\tP$, by (a)),
and in either case $Gn\sim\tN$.

(d) Suppose $G1n\sim\tO$.
Then $G1$, $G11$, and $Gn$
are all $\tN$-games.
$Gn$ must have a $\tP$-option,
but $G1\sim\tN$
and no option $G'n$ or $Gm$ ($2 \leq m < n$)
can be a $\tP$-game
(by (a)),
so $G$ itself must be a $\tP$-game.
Also, since $G11\sim\tN$ and $G1\not\sim\tP$,
there must exist $G'$ with $G'11\sim\tP$.
Then $G'1\sim\tO$ and $G'\sim\tN$,
which is inconsistent with $G\sim\tP$.

(e) Every option of $Gmn$ has a component heap of size 2 or more,
so $G+m+n$ has no $\tP$-options, by (a).

(f) Suppose $Gmn \sim \tO$.
Then $G$ can't be $0$ (by Claim \ref{c:gangup}),
so it must have an option $G'$;
$G' mn \sim \tN$, contradicting (e).

(g) Suppose $Gmn\sim\tP$.
Then $G$ can't be $0$ (by Claim \ref{c:gangup}),
so it must have an option $G'$;
$G' mn \sim \tO$, contradicting (f).
\end{proof}

Note that (e), (f), and (g) together imply
that $Gmn\sim\tQ$ for all $m,n \geq 2$.

\begin{claim} The following are true for all games $G$:
\begin{enumerate}[(A)]
\item $m \equiv n$ for $m,n \geq 3$.
\item $1m \equiv 1n$ for $m,n \geq 2$.
\item $Gmn \equiv 22$ for $m,n \geq 2$.
\end{enumerate}
\end{claim}

\begin{proof}\ \par

(A) Take an arbitrary game $X$.
We know that each of $Xm$, $Xn$
is either of type $\tN$ or type $\tQ$
(by (a), (b) above).
If either of them is a $\tN$-game,
then so is the other (by (c)),
and if neither of them is a $\tN$-game,
then both are $\tQ$-games.
Either way, $m+X$ and $n+X$
have the same type.

(B) The proof is similar, except that one needs
(d) instead of (b).

(C) For all $X$, $Gmn+X=(GX)mn\sim\tQ$
and $22+X = (X)22 \sim \tQ$.
\end{proof}

To reduce a given Nim-position 
$G=n_1+n_2+...+n_r$
to one of the previously tabulated forms,
first replace every $n_i > 3$ by 3.
This puts $G$ in the form
$1^a 2^b 3^c$.
If $b+c \geq 2$,
then we have $G \equiv 22$.
Otherwise, we have $G$ in the form
$1^a$, $1^a 2$, or $1^a 3$.
Since $13 \equiv 12$,
the last of these cases can be reduced to $1^a 2$
unless $a=0$.

\section{Equivalence Classes}

The Nim game $22$ has the property that
if one adds to it any other Nim-position,
one gets a game of type $\tQ$.
In fact, if one adds {\it any\/} game whatsoever to $22$,
one still gets a game of type $\tQ$.
$22$ is thus an element of an important equivalence class,
consisting of all games $G$ such that
$G+X \sim \tQ$ for all games $X$.
We call this class the {\it equivalence class of infinity\/}.
This equivalence class is a sort of a black hole,
metaphorically speaking;
add any game to the black hole,
and all you get is the black hole.

If you take a two-player game
for which a nice theory exists
and study the three-player version,
then it is unfortunately nearly always the case
that most of the positions in the game
are in the equivalence class of infinity.

There are some games which are ``close'' to infinity.
Paradoxically, such games can give us interesting information
about games that are very far away from infinity.
Consider, for instance, the $\tO$-game $\near = \{2\}$
(the game whose sole option is a Nim-heap of size 2).

\begin{claim}
The only game $G$ for which $G+\near \not\sim \tQ$
is the game 0.
\end{claim}

\begin{proof}
Let $X$ be the simplest game not identical to 0
such that $X+\near \not\sim \tQ$.

Case I: $X+\near \sim \tP$.
Then $X+\far \sim \tO$.
But Claim \ref{c:inequiv}(b),
together with the fact that
$\far$ is equivalent to every Nim-position $n$ with $n \geq 3$,
tells us that this can't happen.

Case II: $X+\near \sim \tN$.
The winning option of $X+\near$ can't be $X+\far$,
by Claim \ref{c:inequiv}(a),
so it must be an option of the form $X'+\near$.
But then $X'+\near \sim \tP$,
which contradicts the assumed minimality of $X$.
($X'=0$ won't help us,
since $0+\near \sim \tO$, not $\tP$.)

Case III: $X+\near \sim \tO$.
Letting $X'$ be any option of $X$,
we have $X'+\near \sim \tN$.
This contradicts the assumed minimality of $X$.
\end{proof}

This implies that no game is equivalent to 0.

\section{Open Questions}

\begin{question} How do the doubling and tripling tables
(Tables \ref{t:doubling} and \ref{t:trebling})
extend to higher compound sums of a game with itself?
\end{question}

\begin{question} Is there a decision procedure for determining when two 
impartial three-player games are equivalent to each other?
\end{question}

\begin{question} What does the ``neighborhood of infinity'' look like?
The game $\near \sim \tO$ has the property that when you
add it to any non-trivial game, you get $\tQ$.
Is there a game of type $\tN$ with this property?
Is there one of type $\tP$ with this property?
\end{question}

\begin{question} How does the theory generalize to $n$ players, with $n > 3$?
It is not hard to show that the portion of Table \ref{t:subtraction} in the upper
left corner generalizes to the case of more than three players
in a straightforward way.  However, carrying the theory beyond
this point seems like a large job.  Here are two particular
questions that seem especially interesting: Can an $n$-fold
sum of a game with itself be a win for any of the players
other than the $n$th?  Does there exist a ``black hole'' $X$ 
such that for all games $Y$, $X+Y$ is a win for any coalition 
with over half the players?
\end{question}

\section{Acknowledgments}

\begin{small}
This research was supported by a Knox Fellowship from Harvard College.
I express deep appreciation to John Conway
for his encouragement and for stimulating conversations.
I also thank Richard Guy and Phil Straffin 
for many helpful remarks on the manuscript.
\end{small}

\newcommand{\BCG}{1}
\newcommand{\Con}{2}
\newcommand{\Li}{3}
\newcommand{\Str}{4}

\end{document}